\documentclass[a4paper,11pt, reqno]{amsart}
\usepackage{amssymb,amsthm,amsmath}
\usepackage{cite}

\usepackage{amsmath,amsthm,amscd,amssymb}
\usepackage{latexsym}
\usepackage[colorlinks,citecolor=red,pagebackref,hypertexnames=false]{hyperref}

\theoremstyle{plain}
\newtheorem{theorem}{Theorem}[section]
\newtheorem{Def}[theorem]{Definition}

\newtheorem{proposition}[theorem]{Proposition}

\theoremstyle{definition}

\theoremstyle{remark}
\newtheorem{remark}[theorem]{Remark}

\newtheorem{case[theorem]}{Case}

\def \R{{\mathbb R}}
\def \N{{\mathbb N}}

\def \Z{{\mathbb Z}}

\def \T{{\mathbb T}}

\def\norm#1.#2.{\lVert#1\rVert_{#2}}

\def\R{\mathbb R}

\title{Localization operators on discrete Orlicz modulation spaces}

\author{Aparajita Dasgupta}
\address{Department of Mathematics, Indian Institute of Technology Delhi, New Delhi 110016, India} 

\email{adasgupta@maths.iitd.ac.in}

\author{Anirudha Poria}
\address{Department of Applied Mathematics, School of Mathematics and Physics, Xi’an Jiaotong-Liverpool University, Suzhou 215123, China}

\email{Anirudha.Poria@xjtlu.edu.cn}



\date{\today}
\begin{document}
\maketitle

\begin{abstract}
In this paper, we introduce Orlicz spaces on $ \mathbb Z^n \times \mathbb T^n $ and Orlicz modulation spaces on $\mathbb Z^n$, and study inclusion relations, convolution relations, and duality of these spaces. We show that the Orlicz modulation space $M^{\Phi}(\mathbb Z^n)$ is close to the modulation space $M^{2}(\mathbb Z^n)$ for some particular Young function $\Phi$. Then, we study localization operators on $\mathbb Z^n$. In particular, using appropriate classes for symbols, we prove that these operators are bounded on Orlicz modulation spaces on $\mathbb Z^n$, compact and in the Schatten--von Neumann classes.
\end{abstract}

%%% ----------------------------------------------------------------------
\maketitle
%%% ----------------------------------------------------------------------
%\tableofcontents

\section{Introduction}  

Operators that localize in time and frequency serve as an important mathematical instrument for examining functions across different areas on the time-frequency plane. These can be considered transformations that alter a function’s characteristics in both time and frequency domains, resulting in a reconstructed filtered signal. Daubechies in \cite{dau88, dau90, pau88}, Ramanathan and Topiwala in \cite{ram93} introduced the time-frequency localization operators, and these operators were extensively investigated in \cite{fei02, won99, won02}. This category of operators is found across diverse fields of both applied and pure mathematics and has attracted the attention of numerous researchers. Recognized as a significant novel mathematical instrument, localization operators have been widely applied in areas such as differential equations theory, signal processing, time-frequency analysis, and quantum mechanics (see \cite{cor03, mar02, now02, gro01, ram93, won02}). These operators are also referred to as Gabor multipliers, anti-Wick operators, Toeplitz operators, or wave packets (see \cite{now02, ber71, cor78, fei02}). For an in-depth exploration of localization operators theory, we direct readers to the series of papers authored by Wong \cite{bog04, liu07, won2001, won2002, won2003}, and the book of Wong \cite{won02}. In this paper, our aim is to explore the localization operators on Orlicz modulation spaces on $\Z^n$.

Localization operators have been characterized through the Schr\"odinger representation and the short-time Fourier transform, indicating their study as components of time-frequency analysis. To gain a deeper insight into these operators, modulation spaces serve as suitable function spaces, given their connection with the short-time Fourier transform. Introduced by Feichtinger in \cite{fei03}, modulation spaces constitute a family of spaces for functions and distributions. Since then, the theory of these spaces has been expanded in various ways (see \cite{gro01}). The concept of modulation spaces was extended and investigated using Orlicz spaces and mixed-norm Orlicz spaces in \cite{sch13}. Orlicz spaces are important types of Banach function spaces that are considered in mathematical analysis. These spaces naturally generalize $L^p$-spaces and contain certain Sobolev spaces as subspaces. Orlicz spaces appear in various computations such as the Zygmund space $L \log^+ L$, which is a Banach space related to Hardy--Littlewood maximal functions. Like many other function spaces, there has been a recent interest in the case of Orlicz modulation spaces. Such spaces are obtained by imposing Orlicz norm estimates on the short-time Fourier transforms of the involved functions and distributions. Since the family of Orlicz spaces contains all Lebesgue spaces, the family of Orlicz modulation spaces contain all classical modulation spaces. In particular, the Orlicz modulation spaces are a subfamily of broader classes of modulation spaces. For a detailed study on comparisons of the existing modulation spaces and Orlicz modulation spaces, we refer to \cite{gum24}. Some recent investigations on Orlicz modulation spaces can be found in \cite{sch13, tof22}. In this paper, we introduce Orlicz spaces on $\Z^n \times \T^n$, and Orlicz modulation spaces on $\Z^n$, and study inclusion relations, convolution relations, and duality properties of these spaces. Moreover, we prove that the Orlicz modulation space $M^{\Phi}(\Z^n)$ is close to the modulation space $M^{2}(\Z^n)$ for some particular Young function $\Phi$.

Given that localization operators fall within the category of pseudo-differential operators, recent research on pseudo-differential operators on $\Z^n$ (see \cite{ bot20}) and $\hbar \Z^n$ (see \cite{bot24}) has inspired us to investigate the localization operator on $\Z^n$. Also, our recent works on localization operators on discrete modulation spaces established a strong connection between the theory of pseudo-differential operators and localization operators on $\Z^n$ (see \cite{das23}). In this paper, our main aim is to introduce the Orlicz modulation space on $\Z^n$ using the short-time Fourier transform on $\Z^n \times \T^n$. Then, using suitable conditions on symbols and windows, we show that the localization operators are bounded on Orlicz modulation spaces on $\mathbb Z^n$, compact and in the Schatten--von Neumann class.

In the exploration of the localization operator on $\mathbb{Z}^n$, an important question emerges regarding the appropriate types of spaces to be used for the symbol. For the localization operator on $\mathbb{R}^{n}$, the symbol is represented by a function on $\mathbb{R}^{n} \times \mathbb{R}^{n} $. Recent studies on pseudo-differential operators on topological groups $G$ indicate that the ideal phase space for operation is $G  \times \widehat{G}$, with $\widehat{G}$ being the dual group of $G$ (see \cite{mol19}). Given that the dual group of $\mathbb{R}^{n}$ is identical to $\mathbb{R}^{n}$ itself, the phase space for defining symbols is thus $\mathbb{R}^{n} \times \mathbb{R}^{n}$. In the case of the group $\mathbb{Z}^n$, its dual group is $\T^n$,  making the phase space $\mathbb{Z}^n \times \T^n$. In our paper, we consider the symbol as a function defined on $\mathbb{Z}^n \times \T^n$ and focus our investigation on the localization operator on $\mathbb{Z}^n$.

This paper is organized as follows. In Section \ref{sec2}, we recall some basic facts on the short-time Fourier transform and Orlicz spaces. Also, we study the mixed Orlicz spaces on $ \Z^n \times \T^n $ parameterized with two (quasi-)Young functions. In Section \ref{sec3}, we introduce the Orlicz modulation spaces on $\Z^n$, and study inclusion relations, convolution relations, and duality properties of these spaces. Then, we show that the Orlicz modulation space $M^{\Phi}(\Z^n)$ is close to the modulation space $M^{2}(\Z^n)$ for some particular Young function $\Phi$. Finally, in Section \ref{sec4}, using suitable modulation spaces, Lebesgue spaces or mixed Orlicz spaces on $\Z^n \times \T^n$ as appropriate classes for symbols, we study the localization operator on $\Z^n$ and show that these operators are bounded on Orlicz modulation spaces on $\Z^n$, compact and in the Schatten--von Neumann class.
 
\section{Preliminaries}\label{sec2} 

In this section, we recall some necessary definitions and results related to the short-time Fourier transform and Orlicz spaces. Let us start with some basic definitions.

\begin{Def}
Let $1 \leq p<\infty$. 
\begin{enumerate}
\item We define $\ell^{p}(\mathbb{Z}^n)$ to be the set of all measurable functions $F$ on $\mathbb{Z}^n$ such that  
\[ \|F\|_{\ell^{p}(\mathbb{Z}^n)}^{p}=\sum_{k \in \mathbb{Z}^n}|F(k)|^{p}<\infty. \]

\item We define $L^{p}(\T^n)$ to be the set of all measurable functions $f$ on $\T^n$ for which
\[ \|f\|_{L^{p} (\T^n )}^{p}= \int_{\T^n} |f(w)|^{p} \;dw < \infty. \]

\item The Fourier transform $\mathcal{F}_{\mathbb{Z}^n} F$ of $F \in \ell^{1}(\mathbb{Z}^n)$ is the function on $\T^n$, and defined by
\[ \left(\mathcal{F}_{\mathbb{Z}^n} F \right)(w)=\sum_{k \in \mathbb{Z}^n} e^{-2 \pi i k \cdot w} F(k), \quad w \in \T^n. \]

\item Let $f$ be a function on $\T^n$. The Fourier transform $\mathcal{F}_{\T^n} f$ of $f$ is the function on $\mathbb{Z}^n$, and defined by
\[ \left(\mathcal{F}_{\T^n} f\right)(k)= \int_{\T^n} e^{2 \pi i k \cdot w} f(w) \;dw, \quad k \in \mathbb{Z}^n. \]

\item The space of all measurable functions $H$ on $\mathbb{Z}^n \times \T^n$ such that
\[ \|H\|_{L^{p}\left(\mathbb{Z}^n \times \T^n \right)}^{p}= \sum_{k \in \mathbb{Z}^n} \int_{\T^n}  |H(k, w)|^{p} \;dw < \infty \]
is denoted by $L^{p}\left(\mathbb{Z}^n \times \T^n \right)$. 
\end{enumerate}
\end{Def}
Notice that, $\mathcal{F}_{\mathbb{Z}^n}=\mathcal{F}_{\T^n}^{-1}=\mathcal{F}_{\T^n}^{*} $, and for $F \in \ell^{2}(\mathbb{Z}^n)$, $\left\|\mathcal{F}_{\mathbb{Z}^n} F\right\|_{L^{2}\left(\T^n\right)}=\|F\|_{\ell^{2}(\mathbb{Z}^n)}$. Moreover, $\mathcal{F}_{\mathbb{Z}^n}: \ell^{2}(\mathbb{Z}^n) \rightarrow L^{2} (\T^n)$ is a surjective isomorphism.

Next, we define the Schwartz space $\mathcal{S}\left(\mathbb{Z}^{n}\right)$ on $\mathbb{Z}^{n}$ to be the space of rapidly decreasing functions $g: \mathbb{Z}^{n} \rightarrow \mathbb{C}$. That is, $g \in \mathcal{S}\left(\mathbb{Z}^{n}\right)$ if for any $M<\infty$ there exists a constant $C_{g, M}$ for which $$|g(k)| \leq C_{g, M} (1+|k|)^{-M}, \quad \text { for all } k \in \mathbb{Z}^{n}.$$ 
For $j \in \mathbb{N}_{0}= \N \cup \{0\}$, we define the seminorms $p_{j}(g):=\sup\limits_{k \in \mathbb{Z}^{n}}(1+|k|)^{j}|g(k)|$. Then, the topology on $\mathcal{S}\left(\mathbb{Z}^{n}\right)$ is given by seminorms $p_{j}$. Also, we define the space of tempered distributions $\mathcal{S}^{\prime}\left(\mathbb{Z}^{n}\right)$ to be the space of all continuous linear functionals on $\mathcal{S}\left(\mathbb{Z}^{n}\right)$. 
 
Now, we define the short-time Fourier transform (STFT) on $\Z^n \times \T^n$. Let $f \in \ell^2(\Z^n)$, and fix $k \in \Z^n$, $w \in \T^n$. For $m \in \Z^n$, we define the translation operator $T_k$ by $T_k f(m)=f(m-k)$ and the modulation operator $M_w$ by $M_w f(m)=e^{2 \pi i w \cdot m} f(m)$. For a fixed window function $g \in \mathcal{S}\left(\mathbb{Z}^{n}\right)$, we define the STFT of a function $f \in \mathcal{S}^{\prime}\left(\mathbb{Z}^{n}\right)$ with respect to $g$ to be the function on $\Z^n \times \T^n$ given by 
\begin{eqnarray*}
V_g f (m, w)= \left\langle f, M_w T_m g \right\rangle 
=  \sum_{k \in \Z^n} f(k) \overline{M_w T_m g (k)} 
=  \sum_{k \in \Z^n} f(k) \overline{g(k-m)} e^{-2 \pi i w \cdot k}.
\end{eqnarray*} 
Let $\tilde{g}(k)=g(-k)$, for $k \in \Z^n$. Using the convolution on $\Z^n$, we write $V_g f$ as
\[ V_g f (m, w)= e^{-2 \pi i w \cdot m} \left( f * M_w \overline{\tilde{g}} \right)(m). \]
The STFT on $\Z^n \times \T^n$ satisfies the following properties (see \cite{das23}). 
\begin{proposition}
\begin{itemize}
\item[(1)] For any $f_1, f_2, g_1, g_2 \in \ell^2(\Z^n)$,
\begin{equation}\label{eq18}
\left\langle V_{g_1} f_1, V_{g_2} f_2 \right\rangle_{L^2(\Z^n \times \T^n)}=\langle f_1, f_2 \rangle_{\ell^2(\Z^n)} \; \langle g_2, g_1 \rangle_{\ell^2(\Z^n)}.
\end{equation}

\item[(2)] Let $g \in \ell^2(\Z^n)$. For any $f \in \ell^2(\Z^n)$, we have
\begin{equation}\label{eq17}
\left\Vert V_g f \right\Vert_{L^2(\Z^n \times \T^n)} = \Vert f \Vert_{\ell^2(\Z^n)} \; \Vert g \Vert_{\ell^2(\Z^n)}.
\end{equation}

\item[(3)] Let $g, h \in \ell^2(\Z^n)$ and $\langle g, h \rangle_{\ell^2(\Z^n)} \neq 0$. For any $f \in \ell^2(\Z^n)$, we have
\[f=\frac{1}{\langle h, g \rangle_{\ell^2(\Z^n)}}
\sum_{m \in \Z^n} \int_{\T^n} V_g f(m,w) \; M_w T_m h  \; dw.\]
\end{itemize}
\end{proposition}

Next, we define the Orlicz spaces on $\Z^n \times \T^n$. We first need to define the convex function. A function $\Phi:[0, \infty] \rightarrow [0, \infty]$ is called convex if
$$ \Phi\left(x_1 y_1+x_2 y_2\right) \leq x_1 \Phi\left(y_1\right)+x_2 \Phi\left(y_2\right), $$
where $x_j, y_j \in \mathbb{R}$ satisfy $x_j, y_j \geq 0$ for $j=1,2$ and $x_1+x_2=1$. 

We recall the definition of the Young function and quasi-Young function (see \cite{gum24}).
\begin{Def}
\begin{itemize}
\item[(1)] A function $\Phi : [0, \infty] \rightarrow [0, \infty]$ is called a Young function if $\Phi$ is convex, $\Phi(0)=0$ and $\lim\limits_{t \rightarrow \infty} \Phi(t)=\Phi(\infty)=\infty$.

\item[(2)] A function $\Phi_0 : [0, \infty] \to [0, \infty]$ is called a quasi-Young function of order $p \in (0,1]$ if there is a Young function $\Phi$ such that $\Phi_0(t)=\Phi\left(t^{p}\right)$, where $t \in [0, \infty]$.
\end{itemize}
\end{Def}
Now, we recall the definition of the continuous and discrete Orlicz spaces (see \cite{gum24, mal00}).
\begin{Def}
Let $\Phi$ be a (quasi-)Young function.
\begin{itemize}
\item[(1)] The continuous Orlicz space $L^{\Phi} (\mathbb{R}^n)$ consists of all measurable functions $f: \mathbb{R}^n \rightarrow \mathbb{C}$ such that
$$ \|f\|_{L^{\Phi}\left(\mathbb{R}^n\right)}:=\inf \left\{b>0: \int_{\mathbb{R}^n} \Phi\left(\frac{|f(x)|}{b}\right) d x \leq 1\right\} <\infty .$$

\item[(2)] The discrete Orlicz space $\ell^{\Phi} (\mathbb{Z}^n)$ consists of all measurable functions $F: \mathbb{Z}^n \rightarrow \mathbb{C}$ such that  
$$ \|F\|_{\ell^{\Phi}\left(\mathbb{Z}^n\right)}:=\inf \left\{b>0: \sum_{k \in \mathbb{Z}^n} \Phi\left(\frac{\left|F(k)\right|}{b}\right) \leq 1\right\}<\infty .$$
\end{itemize}
\end{Def}

Note that, if $\Phi(t):=t^p$ for some $p \geq 1$, then $L^{\Phi} (\mathbb{R}^n)=L^p (\mathbb{R}^n)$, the Lebesgue spaces of $p$ the integrable functions on $\mathbb{R}^n$. Hence, the continuous Orlicz spaces are the generalization of the Lebesgue spaces. Also, note that, if $\Phi(t)=t^p$ for some $p \geq 1$, then we get $\ell^{\Phi}(\mathbb{Z}^n)=\ell^p (\mathbb{Z}^n)$. Therefore, the discrete Orlicz spaces are the generalization of $\ell^p (\mathbb{Z}^n)$ spaces. Next, we will define Orlicz spaces on $\Z^n \times \T^n$ parameterized with two (quasi-)Young functions. 

\begin{Def}
Let $\Phi_1$ and $\Phi_2$ be two (quasi-)Young functions. 
\begin{itemize}
\item[(1)] The mixed Orlicz space $L^{\Phi_1, \Phi_2} (\Z^n \times \T^n)$ consists of all measurable functions $F: \Z^n \times \T^n \rightarrow \mathbb{C}$ such that $$\|F\|_{L^{\Phi_1, \Phi_2} ( \Z^n \times \T^n)} := \left\|F_{1}\right\|_{L^{\Phi_2}(\T^n)} < \infty,$$
where
$$ F_{1} \left(w\right)=\left\| F\left(\cdot, w \right) \right\|_{\ell^{\Phi_1} (\Z^n)}.$$

\item[(2)]  The mixed Orlicz space $L_{*}^{\Phi_1, \Phi_2} (\Z^n \times \T^n)$ consists of all measurable functions $F: \Z^n \times \T^n \rightarrow \mathbb{C}$ such that $$\| F \|_{L_{*}^{\Phi_1, \Phi_2}(\Z^n \times \T^n)} := \| G \|_{L^{\Phi_2, \Phi_1} (\T^n \times \Z^n)} < \infty, $$ where $$G(w, m)=F(m, w), \quad m \in \Z^n, \; w \in \T^n.$$
\end{itemize}
\end{Def} 
In this paper, we mainly assume that $\Phi$, $\Phi_1$ and $\Phi_2$ above are Young functions. 

\section{Orlicz modulation spaces on $\Z^n$}\label{sec3}

In this section, we define and study Orlicz modulation spaces on $\Z^n$. Modulation spaces were first introduced by Feichtinger in \cite{fei03, fei97}. To define discrete Orlicz modulation spaces it is essential to revisit the concept of modulation spaces defined on $\Z^n$ (see \cite{das23}). 
\begin{Def}
Let $1 \leq p \leq \infty$ and $g \in \mathcal{S}(\Z^n)$. We define the modulation space $M^{p}(\Z^n)$ to be the space of all tempered distributions $f \in \mathcal{S'}(\Z^n)$ for which $V_g f \in L^{p}(\Z^n \times \T^n)$. The norm on $M^{p}(\Z^n)$ is 
\begin{eqnarray*}
\Vert f \Vert_{M^{p}(\Z^n)} 
= \Vert V_g f \Vert_{L^{p}(\Z^n \times \T^n)} 
= \bigg( \sum_{m \in \Z^n} \int_{\T^n} |V_g f(m,w)|^p \; dw \bigg)^{1/p} < \infty,
\end{eqnarray*}
with the usual adjustments if $p$ is infinite.  
\end{Def}

We have the following inclusions
\[ \mathcal{S}(\Z^n) \subset M^1(\Z^n) \subset M^2(\Z^n)=\ell^2(\Z^n) \subset M^\infty(\Z^n) \subset \mathcal{S'}(\Z^n).\]
In particular, $M^p(\Z^n) \hookrightarrow \ell^p(\Z^n)$ for $1 \leq p \leq 2$, and  $\ell^p(\Z^n) \hookrightarrow M^p(\Z^n)$ for $2 \leq p \leq \infty$. Moreover, for $p < \infty$, $(M^{p}(\Z^n))^{'} =M^{p'}(\Z^n)$, where $p'$ is the conjugate exponent of $p$. 
We have similar inclusion relations for modulation spaces on $\Z^n \times \T^n$, which can be derived by employing techniques similar to those used in the study of modulation spaces on locally compact abelian groups (see \cite{bas22, fei03}). For more information on the properties and applications of modulation spaces, we refer to the book by Gr\"ochenig \cite{gro01}.

\begin{Def}
Fix a non-zero window $g \in \mathcal{S}(\Z^n)$, and $0 < p,q \leq \infty$. Let $\Phi$ and $\Psi$ be (quasi-)Young functions. 
\begin{itemize}
\item[(1)] The discrete modulation spaces $M^{p, q} (\Z^n)$ is set of all $f \in \mathcal{S'}(\Z^n)$ such that
\begin{equation}\label{eq002}
\|f\|_{M^{p, q}(\Z^n)} := \left\|V_g f\right\|_{L^{p, q}(\Z^n \times \T^n)}<\infty.
\end{equation} 
The topology of $M^{p, q}(\Z^n)$ is induced by the norm (\ref{eq002}).

\item[(2)] The discrete Orlicz modulation spaces $M^{\Phi} (\Z^n), M^{\Phi, \Psi}(\Z^n)$ and $W^{\Phi, \Psi}(\Z^n)$ are the sets of all $f \in \mathcal{S'}(\Z^n)$ such that
\begin{equation}\label{eq003}
\|f\|_{M^{\Phi}(\Z^n)} := \left\|V_g f\right\|_{L^{\Phi}(\Z^n \times \T^n)}< \infty,  \;\; \|f\|_{M^{\Phi, \Psi}(\Z^n)} := \left\|V_g f\right\|_{L^{\Phi, \Psi}(\Z^n \times \T^n)} < \infty ,
\end{equation} 
 and
\begin{equation}\label{eq004}
\|f\|_{W^{\Phi, \Psi}(\Z^n)} := \left\|V_g f\right\|_{L_{*}^{\Phi, \Psi}(\Z^n \times \T^n)} < \infty,
\end{equation} 
respectively. The topologies of $M^{\Phi}(\Z^n), M^{\Phi, \Psi}(\Z^n)$ and $W^{\Phi, \Psi}(\Z^n)$ are induced by the respective norms in (\ref{eq003}) and (\ref{eq004}).
\end{itemize}
\end{Def}

Note that the definitions of the discrete Orlicz modulation spaces are independent of the choice of the window function $g \in \mathcal{S}(\Z^n)$. In addition, the difference between the spaces $M^{\Phi, \Psi}(\Z^n)$ and $W^{\Phi, \Psi}(\Z^n)$ is in terms of their topological differences. The topology of $M^{\Phi, \Psi}(\Z^n)$ is induced by the mixed Orlicz space $L^{\Phi, \Psi}(\Z^n \times \T^n)$ norm, whereas the topology of $W^{\Phi, \Psi}(\Z^n)$ is induced by the mixed Orlicz space $L_{*}^{\Phi, \Psi}(\Z^n \times \T^n)$ norm.

Next, we present some basic properties of discrete Orlicz modulation spaces. For some recent investigations of Orlicz modulation spaces, we refer to \cite{gum24, tof22}. We use $q$ as the conjugate exponent of $p$ to define the dual space for a Lebesgue space. Similar to the theory of Lebesgue spaces, we can define the complementary function as a counterpart to the conjugate exponent. Also, in the theory of Orlicz spaces, the Young functions are classified using their growth properties. In particular, the $\Delta_2$-condition plays an important role in defining the dual space of an Orlicz space (see \cite{sch13}). 
\begin{Def}
\begin{itemize}
\item[(1)] (Complementary function) Let $\Psi: \R \to \overline{\R^+}$ be defined by $\Psi(y) =\sup\{ x|y|- \Phi(x); x \geq 0\}$. Then $\Psi$ is called the complementary function to the Young function $\Phi$.

\item[(2)] ($\Delta_2$-condition) A Young function $\Phi : \R \to  \R^+$ is said to satisfy the $\Delta_2$-condition, if there exists a constant $C>0$ and $x_0 \in \R^+_0$, such that $\Phi(2x) \leq C \Phi(x)$ for all $x \geq x_0 \geq 0$. The Young function $\Phi$ is said to satisfy local $\Delta_2$-condition, if there are constants $r>0$ and $C>0$ such that $\Phi(2x) \leq C \Phi(x)$ holds when $x \in [0, r]$.
\end{itemize}
\end{Def}
Next, we give a characterisation of the dual space to the Orlicz space on $\Z^n \times \T^n$. If $(\Phi, \Psi)$ is a complementary Young pair and $\Phi$ satisfies a local $\Delta_2$-condition, then $\left( L^{\Phi} (\mathbb{R}^n)\right)^*$ is isometrically isomorphic to $L^{\Psi} (\mathbb{R}^n)$. Similarly, we can show that $\left(\ell^{\Phi} (\Z^n)\right)^*$ is isometrically isomorphic to $\ell^{\Psi} (\Z^n)$. Let $\left(\Phi_i, \Psi_i \right)$ be complementary Young pairs which satisfy local $\Delta_2$-condition and are strictly convex for $i = 1, 2$. Then $\left(L^{\Phi_1, \Phi_2} (\Z^n \times \T^n)\right)^*$ is isometrically isomorphic to  $L^{\Psi_1, \Psi_2} (\Z^n \times \T^n)$. The proofs of these properties can be obtained using a similar method as in \cite{rao91}. 

Let $\left(\Phi_i, \Psi_i \right)$ be complementary Young pairs for $i = 1, 2$. If $f \in \ell^{\Phi_1} (\Z^n)$ and $g \in \ell^{\Psi_1} (\Z^n)$, then we have the following H\"older's inequality for the Orlicz spaces
\begin{equation}\label{eq1}
\| fg \|_{\ell^{1} (\Z^n)} \leq \|f\|_{\ell^{\Phi_1} (\Z^n)} \|g\|_{\ell^{\Psi_1} (\Z^n)}.
\end{equation}
In addition, if we assume that $\Phi_2$ satisfies a local $\Delta_2$-condition, then for $F \in L^{\Phi_1, \Phi_2} (\Z^n \times \T^n) $ and $G \in L^{\Psi_1, \Psi_2} (\Z^n \times \T^n)$, we have the following H\"older's inequality for the Orlicz spaces 
\begin{equation}\label{eq2}
\| FG \|_{L^{1} (\Z^n \times \T^n)} \leq \|F\|_{L^{\Phi_1, \Phi_2} (\Z^n \times \T^n)} \|G\|_{ L^{\Psi_1, \Psi_2} (\Z^n \times \T^n)}.
\end{equation} 
The proofs of inequalities (\ref{eq1}) and (\ref{eq2}) can be obtained using a similar method as in \cite{rao91}. If $\Phi$ is continuous, then the Schwartz class $\mathcal{S}(\Z^n)$ is embedded into the Orlicz space $\ell^{\Phi} (\Z^n)$. Also, if the complementary function $\Psi$ is continuous then the functions in the Orlicz space define tempered distributions on $\Z^n$. More precisely, let $\left(\Phi_i, \Psi_i \right)$ be complementary Young pairs and $\Phi_i$ be continuous for $i = 1, 2$, then we have the following inclusions 
$$\mathcal{S}(\Z^n) \subset \ell^{\Phi_1} (\Z^n) \subset \mathcal{S'}(\Z^n), $$
if $\Psi_1 $ is continuous. Also,
$$\mathcal{S}(\Z^n \times \T^n) \subset L^{\Phi_1, \Phi_2} (\Z^n \times \T^n) \subset \mathcal{S'}(\Z^n \times \T^n), $$ 
if $\Psi_1, \Psi_2 $ are continuous. If $\Phi_1$ and $\Phi_2$ are (quasi-)Young functions, then $L^{\Phi_1, \Phi_2} (\Z^n \times \T^n)$ is translation invariant, which leads to the fact that $M^{\Phi_1, \Phi_2}(\Z^n)$ is translation and modulation invariant (see \cite{tof22}). Using the fact that $\left(\ell^{\Phi} (\Z^n)\right)^*$ is isometrically isomorphic to $\ell^{\Psi} (\Z^n)$, we can extend the convolution relation $\ell^{1} (\Z^n) * \ell^{p} (\Z^n) \subset \ell^{p} (\Z^n)$ to the Orlicz spaces. Here, we present the following convolution relations.
\begin{proposition}\label{pro6}
\begin{enumerate}
\item If $F \in L^1(\Z^n \times \T^n)$, $G \in L^{\Phi_1, \Phi_2} (\Z^n \times \T^n)$, $\Phi_i$ satisfy local $\Delta_2$-condition and strictly convex Young functions for $i = 1, 2$, then
\begin{equation}\label{eq02}
\|F * G\|_{L^{\Phi_1, \Phi_2} (\Z^n \times \T^n)} \leq \|F\|_{L^1(\Z^n \times \T^n)} \; \|G\|_{L^{\Phi_1, \Phi_2} (\Z^n \times \T^n)}.
\end{equation}

\item If $F \in L^{1}(\Z^n \times \T^n)$,  $G \in L^{\Phi} (\Z^n \times \T^n)$ and $\Phi$ satisfies a local $\Delta_2$-condition, then
\begin{equation}\label{eq01}
\|F * G\|_{L^{\Phi} (\Z^n \times \T^n)} \leq \|F\|_{L^1(\Z^n \times \T^n)} \; \|G\|_{L^{\Phi} (\Z^n \times \T^n)}.
\end{equation}
\end{enumerate}
\end{proposition}
\begin{proof}
(1) Let $\left(\Phi_i, \Psi_i \right)$ be complementary Young pairs which satisfy local $\Delta_2$-condition and are strictly convex for $i = 1, 2$. Then $\left(L^{\Phi_1, \Phi_2} (\Z^n \times \T^n)\right)^*$ is isometrically isomorphic to  $L^{\Psi_1, \Psi_2} (\Z^n \times \T^n)$. If $G \in L^{\Phi_1, \Phi_2} (\Z^n \times \T^n)$, then $T_{(l,x)} G \in L^{\Phi_1, \Phi_2} (\Z^n \times \T^n)$ and $\| T_{(l,x)} G \|_{L^{\Phi_1, \Phi_2} (\Z^n \times \T^n)}= \| G \|_{L^{\Phi_1, \Phi_2} (\Z^n \times \T^n)}$. Let $H \in  L^{\Psi_1, \Psi_2} (\Z^n \times \T^n)$. Using H\"older's inequality (\ref{eq2}), we obtain
\begin{align*}
|\langle F * G, H \rangle| 
& = \left| \sum_{m \in \Z^n}  \int_{\T^n} F*G(m,w) \; \overline{H(m,w)} \; dw \right| \\
& \leq \sum_{m \in \Z^n}  \int_{\T^n} \left( \sum_{l \in \Z^n}  \int_{\T^n} \left| G(m-l, w-x) \right|\; |F(l,x)| \; dx \right) |H(m,w)| \; dw \\
& =  \sum_{l \in \Z^n}  \int_{\T^n} \left( \sum_{m \in \Z^n}  \int_{\T^n} \left| T_{(l,x)} G(m, w) \right|\; |H(m,w)| \; dw \right) |F(l,x)| \; dx  \\
& \leq \sum_{l \in \Z^n}  \int_{\T^n}  |F(l,x)| \; \| T_{(l,x)} G \|_{L^{\Phi_1, \Phi_2} (\Z^n \times \T^n)} \; \|H\|_{L^{\Psi_1, \Psi_2} (\Z^n \times \T^n)} \; dx \\
& = \| G \|_{L^{\Phi_1, \Phi_2} (\Z^n \times \T^n)} \; \|H\|_{L^{\Psi_1, \Psi_2} (\Z^n \times \T^n)} \sum_{l \in \Z^n} \int_{\T^n} |F(l,x)| \; dx \\
& = \|F\|_{L^1(\Z^n \times \T^n)} \; \| G \|_{L^{\Phi_1, \Phi_2} (\Z^n \times \T^n)} \; \|H\|_{L^{\Psi_1, \Psi_2} (\Z^n \times \T^n)}.
\end{align*}
By duality, we get
\begin{align*}
\|F * G\|_{L^{\Phi_1, \Phi_2} (\Z^n \times \T^n)} 
& = \sup\left\{|\langle F * G, H \rangle|:\|H\|_{L^{\Psi_1, \Psi_2} (\Z^n \times \T^n)}  \leq 1 \right\} \\
& \leq \|F\|_{L^1(\Z^n \times \T^n)} \; \|G\|_{L^{\Phi_1, \Phi_2} (\Z^n \times \T^n)}.
\end{align*}
(2) The proof follows similarly as in the first part of the proof by choosing $\Phi_1 = \Phi_2 = \Phi$.
\end{proof}
Now, we study a few properties of the discrete Orlicz modulation spaces. Note that, the definitions of these spaces are independent of the choice of the window function $g$. Also, if the Young function satisfies a local $\Delta_2$-condition, these spaces are Banach spaces. Moreover, if the Young functions are also strictly convex, then the mixed-norm discrete Orlicz modulation spaces are Banach spaces (see \cite[Theorems 6 and 7]{sch13}).
\begin{theorem}\label{th5}
If $\Phi$ satisfies a local $\Delta_2$-condition and its complementary function $\Psi$ is continuous, then $M^{\Phi}(\Z^n)$ is a Banach space. Moreover, if $\left(\Phi_i, \Psi_i \right)$ are complementary Young pairs which satisfy local $\Delta_2$-condition, strictly convex and continuous for $i = 1, 2$, then $M^{\Phi_1, \Phi_2}(\Z^n)$ is a Banach space.
\end{theorem} 
The proof of this theorem can be obtained using a similar method as discussed in \cite{sch13}. So, we skip the proof here. Next, we discuss the duality properties of the discrete Orlicz modulation spaces. If $\left(\Phi, \Psi \right)$ is a complementary Young pair, and $\Phi$ satisfies a local $\Delta_2$-condition and continuous, then $(M^{\Phi}(\Z^n))^* \cong M^{\Psi}(\Z^n)$ under the duality relation
\[ \langle f,h \rangle = \left\langle V_{g_0} f, V_{g_0} h \right\rangle = \sum_{m \in \Z^n} \int_{\T^n} V_{g_0} f(m,w) \; \overline{V_{g_0} h(m,w)} \; dw \] 
for $f \in M^{\Phi}(\Z^n)$ and $h \in M^{\Psi}(\Z^n)$, $g_0 \in \mathcal{S}(\Z^n)$. Note that the duality relation is independent of the choice of the window function $g_0$. If $\left(\Phi_i, \Psi_i \right)$ are complementary Young pairs which satisfy local $\Delta_2$-condition, strictly convex and continuous for $i = 1, 2$, then $(M^{\Phi_1, \Phi_2}(\Z^n))^* \cong M^{\Psi_1, \Psi_2}(\Z^n)$ under the duality relation
\[ \langle f,h \rangle = \left\langle V_{g_0} f, V_{g_0} h \right\rangle = \sum_{m \in \Z^n} \int_{\T^n} V_{g_0} f(m,w) \; \overline{V_{g_0} h(m,w)} \; dw \]   
for $f \in M^{\Phi_1, \Phi_2}(\Z^n)$ and $h \in M^{\Psi_1, \Psi_2}(\Z^n)$, $g_0 \in \mathcal{S}(\Z^n)$.
\begin{theorem}\label{th6}
If $\Phi_i$ and $ \Psi_i$ are (quasi-)Young functions such that 
\[ \lim_{x \to 0^+} \frac{\Psi_i(x)}{\Phi_i(x)} \]
exist and are finite for $i = 1, 2$, then
\begin{equation}\label{eq3}
L^{\Phi_1, \Phi_2} (\Z^n \times \T^n) \hookrightarrow  L^{\Psi_1, \Psi_2} (\Z^n \times \T^n) \quad \text{and} \quad M^{\Phi_1, \Phi_2}(\Z^n) \hookrightarrow M^{\Psi_1, \Psi_2}(\Z^n).
\end{equation}
\end{theorem}
\begin{proof}
The proof of the theorem follows similarly as in \cite[Theorem 5.10]{tof22}.
\end{proof}   
\begin{theorem}\label{th7}
Let $\Phi_i$, $ \Psi_i$, $i = 1, 2$ be (quasi-)Young functions. Then the following conditions are equivalent:
\begin{enumerate}
\item $M^{\Phi_1, \Phi_2}(\Z^n) \subseteq M^{\Psi_1, \Psi_2}(\Z^n)$.

\item $ L^{\Phi_1, \Phi_2} (\Z^n \times \T^n) \subseteq  L^{\Psi_1, \Psi_2} (\Z^n \times \T^n) $.

\item There is a constant $x_0 > 0$ such that $\Psi_i(x) \lesssim \Phi_i(x)$ for all $0 \leq x \leq x_0$.
\end{enumerate}
\end{theorem}
\begin{proof}
Conditions (1) and (3), and (2) and (3) are equivalent follows from Theorem \ref{th6} and \cite[Proposition 5.11]{tof22}. Now, we prove that conditions (1) and (2) are equivalent. Using the definition of discrete Orlicz modulation spaces, we have
\[ \|f\|_{M^{\Phi_1, \Phi_2}(\Z^n)} = \left\|V_g f\right\|_{L^{\Phi_1, \Phi_2}(\Z^n \times \T^n)} \; \text{and}  \; \|f\|_{M^{\Psi_1, \Psi_2}(\Z^n)} = \left\|V_g f\right\|_{L^{\Psi_1, \Psi_2}(\Z^n \times \T^n)}. \]
If $M^{\Phi_1, \Phi_2}(\Z^n) \subseteq M^{\Psi_1, \Psi_2}(\Z^n)$, then
\begin{align*}
\|f\|_{M^{\Psi_1, \Psi_2}(\Z^n)} \leq \|f\|_{M^{\Phi_1, \Phi_2}(\Z^n)}
& \Rightarrow  \left\|V_g f\right\|_{L^{\Psi_1, \Psi_2}(\Z^n \times \T^n)} \leq \left\|V_g f\right\|_{L^{\Phi_1, \Phi_2}(\Z^n \times \T^n)}  \\
& \Rightarrow  L^{\Phi_1, \Phi_2} (\Z^n \times \T^n) \subseteq  L^{\Psi_1, \Psi_2} (\Z^n \times \T^n) .
\end{align*}
Similarly, if $ L^{\Phi_1, \Phi_2} (\Z^n \times \T^n) \subseteq  L^{\Psi_1, \Psi_2} (\Z^n \times \T^n) $, then 
\begin{align*}
\left\|V_g f\right\|_{L^{\Psi_1, \Psi_2}(\Z^n \times \T^n)} \leq \left\|V_g f\right\|_{L^{\Phi_1, \Phi_2}(\Z^n \times \T^n)}  
& \Rightarrow \|f\|_{M^{\Psi_1, \Psi_2}(\Z^n)} \leq \|f\|_{M^{\Phi_1, \Phi_2}(\Z^n)} \\
& \Rightarrow  M^{\Phi_1, \Phi_2}(\Z^n) \subseteq M^{\Psi_1, \Psi_2}(\Z^n).
\end{align*}
This completes the proof.
\end{proof}
Note that the constant of the estimate in condition (3) is uniform in $i$, and the converse implication holds without additional assumptions. Next, we show that $M^{\Phi}(\Z^n)$ is close to $M^{2}(\Z^n)$ in some sense. In the following proposition, we consider the same Young function $\Phi$ as in \cite{gum24}, since this function plays an important role in finding the inclusion relations between the Orlicz modulation space and modulation spaces (see \cite{gum24}).
\begin{proposition}\label{pro7}
Let $\Phi$ be a Young function which satisfies 
\begin{equation}\label{eq5}
\Phi(x)= - x^2 \log x, \quad  0 \leq x \leq e^{-\frac{2}{3}} .
\end{equation} 
Then 
\begin{equation}\label{eq4}
M^{p}(\Z^n) \subseteq M^{\Phi}(\Z^n)  \subseteq M^{2}(\Z^n),  \quad p<2,
\end{equation}
with continuous and dense inclusions.
\end{proposition}
\begin{proof}
Using similar arguments as in \cite[Lemma 3.2]{gum24} and Theorem \ref{th7}, we obtain that the inclusions in (\ref{eq4}) hold and are continuous. Since $M^{p}(\Z^n)$, $ p < 2$ is dense in $M^{2}(\Z^n)$, it also follows that $M^{\Phi}(\Z^n)$ is dense in $M^{2}(\Z^n)$.
\end{proof} 
Throughout the following section, we assume that the Young function $\Phi$ satisfies (\ref{eq5}) so that we can use the inclusion relations in (\ref{eq4}) in the proof of the main results.

\section{Localization operators on Orlicz modulation spaces on $\Z^n$}\label{sec4}

Here, we study the localization operators on $\Z^n$ and prove their boundedness. Furthermore, we demonstrate the compactness of these operators and their inclusion in the Schatten--von Neumann class.

\begin{Def}
Let $\sigma \in L^1(\Z^n \times \T^n)  \cup  L^\infty(\Z^n \times \T^n)$. For the symbol $\sigma$ and two window functions $g_1, g_2 \in \mathcal{S} (\mathbb Z^n)$, the localization operator $\mathfrak{L}^{g_1, g_2}_{\sigma}$ is defined on $\ell^2(\Z^n)$ by 
\begin{equation}\label{eq41}
\mathfrak{L}^{g_1, g_2}_{\sigma}f(k)=\sum_{m \in \Z^n}  \int_{\T^n} \sigma(m, w) \; V_{g_1}f(m,w) \; M_w T_m g_2(k) \; dw, \quad k \in \Z^n. 
\end{equation} 
For any $f,  h \in \ell^2(\Z^n)$, we rewrite the operator $\mathfrak{L}^{g_1, g_2}_{\sigma}$ in a weak sense as 
\begin{equation}\label{eq42}
\left\langle \mathfrak{L}^{g_1, g_2}_{\sigma}f , h \right\rangle_{\ell^2(\Z^n)} = \sum_{m \in \Z^n}  \int_{\T^n} \sigma(m, w) \; V_{g_1}f(m,w) \; \overline{V_{g_2}h(m,w)} \; dw.
\end{equation}
\end{Def}

For $1 \leq p \leq \infty$, we define $\mathcal{B}(\ell^p(\Z^n))$ to be the space of all bounded linear operators from $\ell^p(\Z^n)$ into itself. For $p=2$, the space $\mathcal{B}(\ell^2(\Z^n))$ is the C$^*$-algebra of bounded linear operator $\mathcal{A}$ from $\ell^2(\Z^n)$ into itself, equipped with the norm 
\[ \Vert \mathcal{A} \Vert_{\mathcal{B}(\ell^2(\Z^n))}= \sup_{\Vert f \Vert_{\ell^2(\Z^n)} \leq 1} \Vert \mathcal{A}(f) \Vert_{\ell^2(\Z^n)}. \]
To define the Schatten--von Neumann class $S_p$ on $\Z^n$, we need to first recall the definition of singular values of an operator. For a compact operator $\mathcal{A} \in \mathcal{B}(\ell^2(\Z^n))$, the singular values of $\mathcal{A}$ are the eigenvalues of the positive self-adjoint operator $|\mathcal{A}|=\sqrt{\mathcal{A}^* \mathcal{A}}$ and denoted by $\{ s_n(\mathcal{A}) \}_{n \in \mathbb{N}}$. For $1 \leq p < \infty$, we define the Schatten--von Neumann class $S_p$ to be the space of all compact operators whose singular values lie in $\ell^p$, and equipped with the norm 
\[ \Vert \mathcal{A} \Vert_{S_p}= \left(\sum_{n=1}^\infty  (s_n(\mathcal{A}))^p \right)^{1/p}. \]
For $p=\infty$, $S_\infty$ is the class of all compact operators with the norm  $\Vert \mathcal{A} \Vert_{S_\infty} :=\Vert \mathcal{A} \Vert_{\mathcal{B}(\ell^2(\Z^n))}$. For $p = 1$, the trace of an operator $\mathcal{A} \in S_1$ is defined by 
\[ tr(\mathcal{A}) = \sum_{n=1}^\infty \langle \mathcal{A}  v_n, v_n  \rangle_{\ell^2(\Z^n)}, \]
where $\{v_n \}_n$ is any orthonormal basis of $\ell^2(\Z^n)$. In addition, if $\mathcal{A}$ is positive, then 
\[ tr(\mathcal{A})=\Vert \mathcal{A} \Vert_{S_1} .\]
For a compact operator $\mathcal{A}$ on the Hilbert space $\ell^2(\Z^n)$ if the positive operator $\mathcal{A}^*  \mathcal{A} \in S_1$, then we call the operator $\mathcal{A}$ as a Hilbert--Schmidt operator. For any orthonormal basis $\{v_n \}_n$ of $\ell^2(\Z^n)$, we have 
\[ \Vert \mathcal{A}\Vert_{HS}^2 :=\Vert \mathcal{A}\Vert_{S_2}^2= \Vert \mathcal{A}^* \mathcal{A} \Vert_{S_1}= tr(\mathcal{A}^* \mathcal{A}) = \sum_{n=1}^\infty \Vert \mathcal{A} v_n\Vert^2_{\ell^2(\Z^n)}. \]

\subsection{Boundedness and compactness of $\mathfrak{L}^{g_1, g_2}_{\sigma}$}

Here, we consider $g_1, g_2 \in  M^{\Phi}(\Z^n)$, and prove the results related to the boundedness and compactness of $\mathfrak{L}^{g_1, g_2}_{\sigma}$.
\begin{proposition}\label{pro01}
Let $\Phi$ be a Young function, $\sigma \in L^\infty(\Z^n \times \T^n)$ and $g_1, g_2 \in  M^{\Phi}(\Z^n)$. Then, the operator $\mathfrak{L}^{g_1, g_2}_{\sigma} \in \mathcal{B}(\ell^2(\Z^n))$, and
\[ \left\Vert \mathfrak{L}^{g_1, g_2}_{\sigma} \right\Vert_{\mathcal{B}(\ell^2(\Z^n))} \leq \Vert \sigma \Vert_{L^\infty(\Z^n \times \T^n)} \; \Vert g_1 \Vert_{M^{\Phi}(\Z^n)} \; \Vert g_2 \Vert_{M^{\Phi}(\Z^n)}. \]
\end{proposition}
\begin{proof}
Let $f, h \in \ell^2(\Z^n)$. Applying H\"older's inequality, we get
\begin{eqnarray*}
\left| \left\langle \mathfrak{L}^{g_1, g_2}_{\sigma}f, h \right\rangle_{\ell^2(\Z^n)} \right| 
& \leq & \sum_{m \in \Z^n}  \int_{\T^n} |\sigma(m, w)| \; \left|V_{g_1}f(m,w) \right|\; \left|V_{g_2}h(m,w) \right| \; dw \\ 
& \leq & \left\Vert \sigma \right\Vert_{L^\infty(\Z^n \times \T^n)} \left\Vert V_{g_1}f  \right\Vert_{L^2(\Z^n \times \T^n)} \left\Vert V_{g_2}h \right\Vert_{L^2(\Z^n \times \T^n)}.
\end{eqnarray*}
Applying Plancherel's formula (\ref{eq17}), we obtain
\begin{eqnarray*}
\left| \left\langle \mathfrak{L}^{g_1, g_2}_{\sigma}f, h \right\rangle_{\ell^2(\Z^n)} \right| 
\leq  \left\Vert \sigma \right\Vert_{L^\infty(\Z^n \times \T^n)} \Vert f \Vert_{\ell^2(\Z^n)} \; \Vert g_1 \Vert_{\ell^2(\Z^n)} \; \Vert h  \Vert_{\ell^2(\Z^n)} \; \Vert g_2 \Vert_{\ell^2(\Z^n)}.
\end{eqnarray*}
Since $M^{\Phi}(\Z^n) \subset \ell^2(\Z^n)$, we get
\[\Vert g_1 \Vert_{\ell^2(\Z^n)} \leq \Vert g_1 \Vert_{M^{\Phi}(\Z^n)} \quad \text{and} \quad  \Vert g_2  \Vert_{\ell^2(\Z^n)} \leq \Vert g_2 \Vert_{M^{\Phi}(\Z^n)}. \]
Therefore,
\[ \left\Vert \mathfrak{L}^{g_1, g_2}_{\sigma} \right\Vert_{\mathcal{B}(\ell^2(\Z^n))} \leq \Vert \sigma \Vert_{L^\infty(\Z^n \times \T^n)}\; \Vert g_1 \Vert_{M^{\Phi}(\Z^n)} \; \Vert g_2 \Vert_{M^{\Phi}(\Z^n)}. \]
\end{proof}

\begin{proposition}\label{pro1}
Let $\Phi$ be a Young function, $\sigma \in M^1(\Z^n \times \T^n)$ and $g_1, g_2 \in  M^{\Phi}(\Z^n)$. Then, the operator $\mathfrak{L}^{g_1, g_2}_{\sigma} \in \mathcal{B}(\ell^2(\Z^n))$, and 
$$\left\Vert \mathfrak{L}^{g_1, g_2}_{\sigma} \right\Vert_{\mathcal{B}(\ell^2(\Z^n))} \leq \Vert \sigma \Vert_{M^1(\Z^n \times \T^n)} \; \Vert g_1 \Vert_{M^{\Phi}(\Z^n)} \; \Vert g_2 \Vert_{M^{\Phi}(\Z^n)}. $$
\end{proposition}
\begin{proof}
Let $f , h \in \ell^2(\Z^n)$. Using the duality between the modulation spaces $M^\infty(\Z^n \times \T^n)$ and $M^1(\Z^n \times \T^n)$, we obtain
\begin{eqnarray}\label{eq45}
\left| \left\langle \mathfrak{L}^{g_1, g_2}_{\sigma}f, h \right\rangle_{\ell^2(\Z^n)} \right|   
& \leq & \sum_{m \in \Z^n}  \int_{\T^n} |\sigma(m, w)| \; \left|V_{g_1}f(m, w) \; \overline{V_{g_2}h(m, w)} \right|  \; dw \nonumber \\
& \leq & \Vert \sigma \Vert_{M^1(\Z^n \times \T^n)} \left\Vert V_{g_1}f \cdot \overline{V_{g_2}h} \right\Vert_{M^\infty(\Z^n \times \T^n)}.
\end{eqnarray}
Since $L^2(\Z^n \times \T^n) \subset M^\infty(\Z^n \times \T^n)$ and $M^{\Phi}(\Z^n) \subset \ell^2(\Z^n)$, applying Plancherel's formula (\ref{eq17}), we get 
\begin{eqnarray}\label{eq46}
&& \left\Vert V_{g_1}f \cdot \overline{V_{g_2}h} \right\Vert_{M^\infty(\Z^n \times \T^n)} \nonumber \\
& \leq & \left\Vert V_{g_1}f \cdot \overline{V_{g_2}h} \right\Vert_{L^2(\Z^n \times \T^n)} \nonumber \\
& \leq & \left\Vert V_{g_1}f \right\Vert_{L^2(\Z^n \times \T^n)} \left\Vert V_{g_2}h \right\Vert_{L^2(\Z^n \times \T^n)} \nonumber \\
& = & \Vert f \Vert_{\ell^2(\Z^n)} \; \Vert g_1 \Vert_{\ell^2(\Z^n)} \; \Vert h \Vert_{\ell^2(\Z^n)} \; \Vert g_2 \Vert_{\ell^2(\Z^n)} \nonumber \\
& \leq & \Vert f \Vert_{\ell^2(\Z^n)} \; \Vert h \Vert_{\ell^2(\Z^n)} \; \Vert g_1 \Vert_{M^{\Phi}(\Z^n)} \; \Vert g_2 \Vert_{M^{\Phi}(\Z^n)}.
\end{eqnarray}
From (\ref{eq45}) and (\ref{eq46}), we have
\[\left\Vert \mathfrak{L}^{g_1, g_2}_{\sigma} \right\Vert_{\mathcal{B}(\ell^2(\Z^n))} \leq \Vert \sigma \Vert_{M^1(\Z^n \times \T^n)} \;  \Vert g_1 \Vert_{M^{\Phi}(\Z^n)} \; \Vert g_2 \Vert_{M^{\Phi}(\Z^n)} .\]
\end{proof}
\begin{proposition}\label{pro2}
Let $\Phi$ be a Young function, $\sigma \in M^2(\Z^n \times \T^n)$ and $g_1, g_2 \in  M^{\Phi}(\Z^n)$. Then, the operator $\mathfrak{L}^{g_1, g_2}_{\sigma} \in \mathcal{B}(\ell^2(\Z^n))$, and 
\[ \left\Vert \mathfrak{L}^{g_1, g_2}_{\sigma} \right\Vert_{\mathcal{B}(\ell^2(\Z^n))} \leq \Vert \sigma \Vert_{M^2(\Z^n \times \T^n)} \; \Vert g_1 \Vert_{M^{\Phi}(\Z^n)} \; \Vert g_2 \Vert_{M^{\Phi}(\Z^n)}. \]
\end{proposition}  
\begin{proof}
Let $f, h \in \ell^2(\Z^n)$. Applying H\"older's inequality, we get
\begin{eqnarray*}
\left| \left\langle \mathfrak{L}^{g_1, g_2}_{\sigma}f, h \right\rangle_{\ell^2(\Z^n)} \right| 
& \leq & \sum_{m \in \Z^n}  \int_{\T^n} |\sigma(m, w)| \; \left|V_{g_1}f(m, w) \; \overline{V_{g_2}h(m, w)} \right| \; dw \\ 
& \leq & \left\Vert \sigma \right\Vert_{L^2(\Z^n \times \T^n)} \left\Vert V_{g_1}f  \cdot \overline{V_{g_2}h} \right\Vert_{L^2(\Z^n \times \T^n)}.
\end{eqnarray*}
Since $L^2(\Z^n \times \T^n)= M^2(\Z^n \times \T^n)$, using (\ref{eq46}), we obtain
\begin{eqnarray*}
&& \left| \left\langle \mathfrak{L}^{g_1, g_2}_{\sigma}f, h \right\rangle_{\ell^2(\Z^n)} \right| \\
& \leq & \left\Vert \sigma \right\Vert_{M^2(\Z^n \times \T^n)} \Vert f \Vert_{\ell^2(\Z^n)} \; \Vert h \Vert_{\ell^2(\Z^n)} \; \Vert g_1 \Vert_{M^{\Phi}(\Z^n)} \; \Vert g_2 \Vert_{M^{\Phi}(\Z^n)}.
\end{eqnarray*}
Therefore,
\[ \left\Vert \mathfrak{L}^{g_1, g_2}_{\sigma} \right\Vert_{\mathcal{B}(\ell^2(\Z^n))} \leq \Vert \sigma \Vert_{M^2(\Z^n \times \T^n)} \; \Vert g_1 \Vert_{M^{\Phi}(\Z^n)} \; \Vert g_2 \Vert_{M^{\Phi}(\Z^n)}. \]
\end{proof}

\begin{theorem}\label{th1}
Let $1 < p < 2$, $\sigma \in M^p(\Z^n \times \T^n)$, $\Phi$ be a Young function, and $g_1, g_2 \in  M^{\Phi}(\Z^n)$. For fixed $\sigma \in M^p(\Z^n \times \T^n)$, the operator $\mathfrak{L}^{g_1, g_2}_{\sigma}$ can be uniquely extended to a bounded linear operator on $\ell^2(\Z^n)$, for which
\[ \left\Vert \mathfrak{L}^{g_1, g_2}_{\sigma} \right\Vert_{\mathcal{B}(\ell^2(\Z^n))} \leq \Vert \sigma \Vert_{M^p(\Z^n \times \T^n)} \; \Vert g_1 \Vert_{M^{\Phi}(\Z^n)} \; \Vert g_2 \Vert_{M^{\Phi}(\Z^n)}. \] 
\end{theorem} 
\begin{proof}
Let $1< p < 2$ and $\sigma \in M^1(\Z^n \times \T^n) \cap  M^2(\Z^n \times \T^n)$. The modulation spaces $M^{p}$ interpolate similar to the corresponding mixed-norm spaces $L^{p}$. From Proposition \ref{pro1}, for $\sigma \in M^1(\Z^n \times \T^n)$, we have 
$$\left\Vert \mathfrak{L}^{g_1, g_2}_{\sigma} \right\Vert_{\mathcal{B}(\ell^2(\Z^n))} \leq \Vert \sigma \Vert_{M^1(\Z^n \times \T^n)} \; \Vert g_1 \Vert_{M^{\Phi}(\Z^n)} \; \Vert g_2 \Vert_{M^{\Phi}(\Z^n)}.$$
Also, from Proposition \ref{pro2}, for $\sigma \in M^2(\Z^n \times \T^n)$, we have
\[ \left\Vert \mathfrak{L}^{g_1, g_2}_{\sigma} \right\Vert_{\mathcal{B}(\ell^2(\Z^n))} \leq \Vert \sigma \Vert_{M^2(\Z^n \times \T^n)} \; \Vert g_1 \Vert_{M^{\Phi}(\Z^n)} \; \Vert g_2 \Vert_{M^{\Phi}(\Z^n)}. \]
Now, for $1< p < 2$, using the Riesz--Thorin interpolation theorem (see \cite{ste56}), we get 
\[ \left\Vert \mathfrak{L}^{g_1, g_2}_{\sigma} \right\Vert_{\mathcal{B}(\ell^2(\Z^n))} \leq \Vert \sigma \Vert_{M^p(\Z^n \times \T^n)} \; \Vert g_1 \Vert_{M^{\Phi}(\Z^n)} \; \Vert g_2 \Vert_{M^{\Phi}(\Z^n)} . \] 
Let $\sigma \in M^p(\Z^n \times \T^n)$ and $\{ \sigma_n \}_{n \geq 1}$ be a sequence of functions in $M^1(\Z^n \times \T^n) \cap M^2(\Z^n \times \T^n)$ such that $\sigma_n \to \sigma$ in $M^p(\Z^n \times \T^n)$ as $ n \to \infty$. Therefore, for any $n, k \in \mathbb{N}$, we get
\[ \left\Vert \mathfrak{L}^{g_1, g_2}_{\sigma_n} - \mathfrak{L}^{g_1, g_2}_{\sigma_k} \right\Vert_{\mathcal{B}(\ell^2(\Z^n))} \leq \Vert \sigma_n - \sigma_k \Vert_{M^p(\Z^n \times \T^n)} \; \Vert g_1 \Vert_{M^{\Phi}(\Z^n)} \; \Vert g_2 \Vert_{M^{\Phi}(\Z^n)} . \]   
Hence, $\{ \mathfrak{L}^{g_1, g_2}_{\sigma_n} \}_{n \geq 1}$ is a Cauchy sequence in $\mathcal{B}(\ell^2(\Z^n))$. Let $ \mathfrak{L}^{g_1, g_2}_{\sigma_n} \to \mathfrak{L}^{g_1, g_2}_{\sigma}$ as $ n \to \infty$. Then, the limit $ \mathfrak{L}^{g_1, g_2}_{\sigma}$ remains unaffected by the selection of the sequence $\{ \sigma_n \}_{n \geq 1}$, and we have
\begin{eqnarray*}
\left\Vert \mathfrak{L}^{g_1, g_2}_{\sigma} \right\Vert_{\mathcal{B}(\ell^2(\Z^n))} 
& = &  \lim_{n \to \infty} \left\Vert \mathfrak{L}^{g_1, g_2}_{\sigma_n} \right\Vert_{\mathcal{B}(\ell^2(\Z^n))} \\
& \leq & \lim_{n \to \infty} \Vert \sigma_n \Vert_{M^p(\Z^n \times \T^n)} \; \Vert g_1 \Vert_{M^{\Phi}(\Z^n)} \; \Vert g_2 \Vert_{M^{\Phi}(\Z^n)} \\
& =& \Vert \sigma \Vert_{M^p(\Z^n \times \T^n)} \; \Vert g_1 \Vert_{M^{\Phi}(\Z^n)} \; \Vert g_2 \Vert_{M^{\Phi}(\Z^n)}.
\end{eqnarray*}
\end{proof}

\begin{theorem}\label{th2}
Let $1 \leq p \leq 2$, $\sigma \in M^p(\Z^n \times \T^n)$, $\Phi$ be a Young function, and $g_1, g_2 \in  M^{\Phi}(\Z^n)$. Then, the operator $\mathfrak{L}^{g_1, g_2}_{\sigma}: \ell^2(\Z^n) \to \ell^2(\Z^n) $ is compact. 
\end{theorem}  
\begin{proof}
Let $\sigma \in M^1(\Z^n \times \T^n)$ and $\{ v_n \}_n$ be an orthonormal basis for $\ell^2(\Z^n)$. Since $M^1(\Z^n \times \T^n) \subset L^1(\Z^n \times \T^n)$, applying Parseval's identity, we get
\begin{eqnarray*}
&& \sum_{n=1}^\infty \left| \left\langle \mathfrak{L}^{g_1, g_2}_{\sigma} v_n, v_n \right\rangle_{\ell^2(\Z^n)} \right| \\
&& \leq \sum_{n=1}^\infty  \sum_{m \in \Z^n}  \int_{\T^n} \left| \sigma (m,w) \right| \; \left| \langle v_n, M_w T_m g_1 \rangle_{\ell^2(\Z^n)} \right| \; \left| \langle M_w T_m g_2 , v_n \rangle_{\ell^2(\Z^n)} \right| \; dw\\
&& = \sum_{m \in \Z^n}  \int_{\T^n} \left| \sigma (m,w) \right| \left( \sum_{n=1}^\infty  \left| \langle v_n, M_w T_m g_1 \rangle_{\ell^2(\Z^n)} \right|\;  \left| \langle M_w T_m g_2 , v_n \rangle_{\ell^2(\Z^n)} \right| \right) dw \\
&& \leq \frac{1}{2} \sum_{m \in \Z^n}  \int_{\T^n} \left| \sigma (m,w) \right| \left( \sum_{n=1}^\infty |\langle v_n, M_w T_m g_1 \rangle_{\ell^2(\Z^n)}|^2 \right. \\
&& \qquad \qquad \left.  + \sum_{n=1}^\infty |\langle M_w T_m g_2, v_n \rangle_{\ell^2(\Z^n)}|^2 \right) \; dw\\
&& = \frac{1}{2} \Vert \sigma \Vert_{L^1(\Z^n \times \T^n)} \; (\Vert g_1 \Vert^2_{\ell^2(\Z^n)} + \Vert g_2 \Vert^2_{\ell^2(\Z^n)}) \\
&& \leq \frac{1}{2} \Vert \sigma \Vert_{M^1(\Z^n \times \T^n)} \; (\Vert g_1 \Vert^2_{M^{\Phi}(\Z^n)} + \Vert g_2 \Vert^2_{M^{\Phi}(\Z^n)}).   
\end{eqnarray*}
Hence, the operator $\mathfrak{L}^{g_1, g_2}_{\sigma} \in S_1$. Next, let $\sigma \in M^p(\Z^n \times \T^n)$. We consider $\{ \sigma_n \}_{n \geq 1}$ in $M^1(\Z^n \times \T^n) \cap M^2(\Z^n \times \T^n)$ such that $\sigma_n \to \sigma$ in $M^p(\Z^n \times \T^n)$ as $n \to \infty$. Then, applying Theorem \ref{th1}, we obtain
\[ \Vert \mathfrak{L}^{g_1, g_2}_{\sigma_n} - \mathfrak{L}^{g_1, g_2}_{\sigma}  \Vert_{\mathcal{B}(\ell^2(\Z^n))} \leq \Vert \sigma_n - \sigma \Vert_{M^p(\Z^n \times \T^n)} \; \Vert g_1 \Vert_{M^{\Phi}(\Z^n)} \; \Vert g_2 \Vert_{M^{\Phi}(\Z^n)}  \to 0, \]
as $n \to \infty$. Therefore, $\mathfrak{L}^{g_1, g_2}_{\sigma_n} \to \mathfrak{L}^{g_1, g_2}_{\sigma} $ in $\mathcal{B}(\ell^2(\Z^n))$ as $n \to \infty$. From the above, we get that $ \{ \mathfrak{L}^{g_1, g_2}_{\sigma_n} \}_{n \geq 1}$ is a sequence of linear operators in $S_1$ and hence compact, so the operator $\mathfrak{L}^{g_1, g_2}_{\sigma}$ is compact. 
\end{proof}
Now, we calculate the adjoint $(\mathfrak{L}^{g_1, g_2}_{\sigma})^*$ of the operator $\mathfrak{L}^{g_1, g_2}_{\sigma}$ on $\ell^2(\Z^n)$, which determined by the relation
\[ \left\langle \mathfrak{L}^{g_1, g_2}_{\sigma}f , h \right\rangle_{\ell^2(\Z^n)}= \left\langle f , (\mathfrak{L}^{g_1, g_2}_{\sigma})^*h \right\rangle_{\ell^2(\Z^n)}. \] 
We obtain
\begin{eqnarray*}
&& \left\langle \mathfrak{L}^{g_1, g_2}_{\sigma}f , h \right\rangle_{\ell^2(\Z^n)} \\
&=& \sum_{m \in \Z^n}  \int_{\T^n} \sigma(m, w) \; V_{g_1}f(m,w) \; \overline{V_{g_2}h(m,w)} \; dw \\
&=& \sum_{m \in \Z^n}  \int_{\T^n} \sigma(m, w) \; \langle f, M_w T_m g_1 \rangle_{\ell^2(\Z^n)}  \; \overline{V_{g_2}h(m,w)} \; dw 
\end{eqnarray*}
\begin{eqnarray*}
&=& \sum_{m \in \Z^n}  \int_{\T^n} \sigma(m, w) \; \langle f, V_{g_2}h(m,w) \; M_w T_m g_1 \rangle_{\ell^2(\Z^n)}  \; dw \\
&=&  \left\langle f, \sum_{m \in \Z^n}  \int_{\T^n} \overline{\sigma(m, w)} \; V_{g_2}h(m,w) \; M_w T_m g_1 \; dw \right\rangle_{\ell^2(\Z^n)} \\
&=& \left\langle f, \mathfrak{L}^{g_2, g_1}_{\overline{\sigma}}h \right\rangle_{\ell^2(\Z^n)}.
\end{eqnarray*}
Hence, we get
\[ (\mathfrak{L}^{g_1, g_2}_{\sigma})^* = \mathfrak{L}^{g_2, g_1}_{\overline{\sigma}}. \]
Therefore, the operator $\mathfrak{L}^{g_1, g_2}_{\sigma}$ is a self-adjoint operator if $g_1=g_2$ and $\sigma$ is a real-valued function.

\subsection{Localization operators in $S_p$}
     
Here, we show that the operator $\mathfrak{L}^{g, g}_{\sigma} \in S_p$ and give an upper bound of $\Vert \mathfrak{L}^{g, g}_{\sigma} \Vert_{S_p}$. 
\begin{proposition}\label{pro3}
Let $\sigma \in  M^1(\Z^n \times \T^n)$, $\Phi$ be a Young function, and $g \in  M^{\Phi}(\Z^n)$. Then, the operator $\mathfrak{L}^{g, g}_{\sigma}: \ell^2(\Z^n) \to \ell^2(\Z^n) $ is in $S_1$ and 
\[ \Vert \mathfrak{L}^{g, g}_{\sigma} \Vert_{S_1} \leq 4 \Vert \sigma \Vert_{M^1(\Z^n \times \T^n)} \;  \Vert g \Vert^2_{M^{\Phi}(\Z^n)}. \]
\end{proposition}
\begin{proof}
Let $\sigma \in  M^1(\Z^n \times \T^n)$. Then, from Theorem \ref{th2}, the operator $\mathfrak{L}^{g, g}_{\sigma} \in S_1$. Let $\sigma \in M^1(\Z^n \times \T^n)$ be non-negative real-valued. Then $(({\mathfrak{L}^{g, g}_{\sigma}})^* \mathfrak{L}^{g, g}_{\sigma})^{1/2}=\mathfrak{L}^{g, g}_{\sigma}$. Let $\{v_n\}_n$ be an orthonormal basis for $\ell^2(\Z^n)$ consisting of eigenvalues of the operator $(({\mathfrak{L}^{g, g}_{\sigma}})^* \mathfrak{L}^{g, g}_{\sigma})^{1/2}: \ell^2(\Z^n) \to \ell^2(\Z^n)$. If $\sigma$ is non-negative real-valued function, then by using a similar method as in the proof of Theorem \ref{th2}, we get the following estimate
\begin{eqnarray}\label{eq43}
\Vert \mathfrak{L}^{g, g}_{\sigma} \Vert_{S_1}
& =& \sum_{n=1}^\infty \left\langle (({\mathfrak{L}^{g, g}_{\sigma}})^* \mathfrak{L}^{g, g}_{\sigma})^{1/2}v_n, v_n \right\rangle_{\ell^2(\Z^n)} \nonumber \\
&=& \sum_{n=1}^\infty \left\langle \mathfrak{L}^{g, g}_{\sigma} v_n, v_n \right\rangle_{\ell^2(\Z^n)}  
\leq  \Vert \sigma \Vert_{M^1(\Z^n \times \T^n)} \; \Vert g \Vert^2_{M^{\Phi}(\Z^n)} .
\end{eqnarray}
Next, let $\sigma \in M^1(\Z^n \times \T^n)$ be an arbitrary real-valued function. We write $\sigma = \sigma_+ - \sigma_-$, where $\sigma_+=\max(\sigma, 0)$ and $\sigma_-=-\min(\sigma, 0)$. Then, applying the relation (\ref{eq43}), we get
\begin{eqnarray}\label{eq44}
\Vert \mathfrak{L}^{g, g}_{\sigma} \Vert_{S_1}
& = & \Vert \mathfrak{L}^{g, g}_{\sigma_+} -\mathfrak{L}^{g, g}_{\sigma_-} \Vert_{S_1} \nonumber \\
& \leq & \Vert \mathfrak{L}^{g, g}_{\sigma_+} \Vert_{S_1} +\Vert \mathfrak{L}^{g, g}_{\sigma_-} \Vert_{S_1} \nonumber \\
& \leq &  \Vert g \Vert^2_{M^{\Phi}(\Z^n)} (\Vert \sigma_+ \Vert_{M^1(\Z^n \times \T^n)} + \Vert \sigma_- \Vert_{M^1(\Z^n \times \T^n)}) \nonumber \\
& \leq & 2  \Vert g \Vert^2_{M^{\Phi}(\Z^n)} \Vert \sigma \Vert_{M^1(\Z^n \times \T^n)}.
\end{eqnarray}
Finally, let $\sigma \in  M^1(\Z^n \times \T^n)$ be a complex-valued function. Then, we write $\sigma=\sigma_1+ i \sigma_2$, where $\sigma_1, \sigma_2$ are the real and imaginary parts of $\sigma$ respectively. Applying the relation (\ref{eq44}), we get
\begin{eqnarray*}
\Vert \mathfrak{L}^{g, g}_{\sigma} \Vert_{S_1}
& = & \Vert \mathfrak{L}^{g, g}_{\sigma_1}+i \; \mathfrak{L}^{g, g}_{\sigma_2} \Vert_{S_1}  \\
& \leq & \Vert \mathfrak{L}^{g, g}_{\sigma_1} \Vert_{S_1} + \Vert \mathfrak{L}^{g, g}_{\sigma_2} \Vert_{S_1} \\
&\leq & 2  \Vert g \Vert^2_{M^{\Phi}(\Z^n)} (\Vert \sigma_1 \Vert_{M^1(\Z^n \times \T^n)} + \Vert \sigma_2 \Vert_{M^1(\Z^n \times \T^n)}) \nonumber \\
& \leq & 4 \Vert \sigma \Vert_{M^1(\Z^n \times \T^n)}  \Vert g \Vert^2_{M^{\Phi}(\Z^n)}  .
\end{eqnarray*}
\end{proof}

\begin{theorem}
Let $\sigma \in  M^1(\Z^n \times \T^n)$, $\Phi$ be a Young function, and $g \in  M^{\Phi}(\Z^n)$. Then, for $1 \leq p \leq \infty$, the operator $\mathfrak{L}^{g, g}_{\sigma}: \ell^2(\Z^n) \to \ell^2(\Z^n)$ is in $S_p$, and 
\[\Vert \mathfrak{L}^{g, g}_{\sigma} \Vert_{S_p} \leq 2^{2/p} \Vert \sigma  \Vert_{M^1(\Z^n \times \T^n)} \;  \Vert g \Vert^{2}_{M^{\Phi}(\Z^n)}.\]
\end{theorem}
\begin{proof}
Using interpolation theorems (see \cite{won02}, Theorems 2.10 and 2.11), Proposition \ref{pro2}, and Proposition \ref{pro3}, we obtain the proof of the theorem.
\end{proof}
In the following, we obtain a lower bound of $\Vert \mathfrak{L}^{g, g}_{\sigma} \Vert_{S_1}$, and improve the constant given in Proposition \ref{pro3}.
\begin{theorem}\label{th4}
Let $\sigma \in M^1(\Z^n \times \T^n)$ be a non-negative real-valued function, $\Phi$ be a Young function, and $g \in  M^{\Phi}(\Z^n)$. Then, the operator $\mathfrak{L}^{g, g}_{\sigma}$ is in $S_1$, and 
\begin{eqnarray*}
\frac{1}{\Vert g \Vert^2_{M^{\Phi}(\Z^n)}}\; \Vert \tilde{\sigma} \Vert_{L^1(\Z^n \times \T^n)}
\leq \Vert \mathfrak{L}^{g, g}_{\sigma} \Vert_{S_1} \leq \Vert g \Vert^2_{M^{\Phi}(\Z^n)} \; \Vert \sigma \Vert_{M^1(\Z^n \times \T^n)},
\end{eqnarray*}
where $\tilde{\sigma}$ is given by $\tilde{\sigma}(m, w)=\left\langle \mathfrak{L}^{g, g}_{\sigma}(M_w T_m g), M_w T_m g \right\rangle_{\ell^2(\Z^n)}$.
\end{theorem}
\begin{proof}
Let $\sigma \in M^1(\Z^n \times \T^n)$. Then, from Proposition \ref{pro3}, the operator $\mathfrak{L}^{g, g}_{\sigma} \in S_1$. If $\{s_n(\mathfrak{L}^{g, g}_{\sigma}) \}_n$ are the positive singular values of $\mathfrak{L}^{g, g}_{\sigma}$, $\{ v_n \}_n$ is an orthonormal basis for the orthogonal complement of the null space of $\mathfrak{L}^{g, g}_{\sigma}$ consisting of eigenvectors of $|\mathfrak{L}^{g, g}_{\sigma}|$ and $\{ u_n \}_n$ is an orthonormal set in $\ell^2(\Z^n)$, then using the canonical form of compact operators (see \cite{won02}, Theorem 2.2), we get
\begin{equation}\label{eq49}
\mathfrak{L}^{g, g}_{\sigma}f =\sum_{n=1}^{\infty} s_n(\mathfrak{L}^{g, g}_{\sigma}) \langle f, v_n \rangle_{\ell^2(\Z^n)} u_n.
\end{equation}
Therefore, we obtain
\[ \sum_{n=1}^{\infty}  \langle \mathfrak{L}^{g, g}_{\sigma} v_n, u_n \rangle_{\ell^2(\Z^n)}   = \sum_{n=1}^{\infty} s_n(\mathfrak{L}^{g, g}_{\sigma})= \Vert \mathfrak{L}^{g, g}_{\sigma} \Vert_{S_1}.  \]
Applying Bessel's inequality and Cauchy--Schwarz's inequality, we obtain
\begin{eqnarray*}
\Vert \mathfrak{L}^{g, g}_{\sigma} \Vert_{S_1} 
& = &  \sum_{n=1}^{\infty}  \langle \mathfrak{L}^{g, g}_{\sigma} v_n, u_n \rangle_{\ell^2(\Z^n)} \\
& = & \sum_{n=1}^{\infty} \sum_{m \in \Z^n}  \int_{\T^n} \sigma(m, w) \; V_{g} v_n (m,w) \; \overline{V_{g} u_n (m,w)} \; dw \\
& \leq & \sum_{m \in \Z^n}  \int_{\T^n} |\sigma(m, w)| \left(\sum_{n=1}^{\infty} |V_{g} v_n (m,w)|^2 \right)^{1/2} \\ 
&& \qquad \qquad   \left(\sum_{n=1}^{\infty} |V_{g} u_n (m,w)|^2 \right)^{1/2} dw \\
& \leq & \Vert \sigma \Vert_{L^1(\Z^n \times \T^n)} \;\Vert g \Vert^2_{\ell^2(\Z^n)}  \\
& \leq & \Vert \sigma \Vert_{M^1(\Z^n \times \T^n)}\; \Vert g \Vert^2_{M^{\Phi}(\Z^n)}.
\end{eqnarray*}
To prove that $\tilde{\sigma} \in L^1(\Z^n \times \T^n)$, applying formula (\ref{eq49}), we get
\begin{eqnarray*}
&& |\tilde{\sigma}(m, w)| \\
& = & \left|\left\langle \mathfrak{L}^{g, g}_{\sigma}( M_w T_m g),  M_w T_m g \right\rangle_{\ell^2(\Z^n)}\right| \\
& = & \left| \sum_{n=1}^{\infty} s_n(\mathfrak{L}^{g, g}_{\sigma}) \left\langle  M_w T_m g, v_n \right\rangle_{\ell^2(\Z^n)} \left\langle u_n,  M_w T_m g \right\rangle_{\ell^2(\Z^n)} \right| \\
& \leq & \frac{1}{2} \sum_{n=1}^{\infty} s_n(\mathfrak{L}^{g, g}_{\sigma}) \left( \left| \left\langle  M_w T_m g, v_n \right\rangle_{\ell^2(\Z^n)}  \right|^2 + \left| \left\langle  M_w T_m g, u_n  \right\rangle_{\ell^2(\Z^n)} \right|^2 \right). 
\end{eqnarray*}  
Now, using Plancherel's formula (\ref{eq17}), we obtain
\begin{eqnarray*}
&& \Vert \tilde{\sigma} \Vert_{L^1(\Z^n \times \T^n)} 
= \sum_{m \in \Z^n}  \int_{\T^n} |\tilde{\sigma}(m, w)| \; dw  \\
&& \leq \frac{1}{2} \sum_{n=1}^{\infty} s_n(\mathfrak{L}^{g, g}_{\sigma}) \sum_{m \in \Z^n} \int_{\T^n} \left(  \left| \left\langle  M_w T_m g, v_n \right\rangle_{\ell^2(\Z^n)} \right|^2  \right. \\
&& \qquad \qquad   + \left. \left| \left\langle  M_w T_m g, u_n  \right\rangle_{\ell^2(\Z^n)} \right|^2  \right) dw  \\
&& \leq  \Vert g \Vert^2_{M^{\Phi}(\Z^n)} \sum_{n=1}^{\infty} s_n(\mathfrak{L}^{g, g}_{\sigma})  \\
&& = \Vert g \Vert^2_{M^{\Phi}(\Z^n)} \; \Vert \mathfrak{L}^{g, g}_{\sigma} \Vert_{S_1} .
\end{eqnarray*}
\end{proof} 

\subsection{$M^{\Phi}(\Z^n)$ Boundedness} 
 
Here, we obtain that the operators $\mathfrak{L}^{g_1, g_2}_{\sigma} : M^{\Phi}(\Z^n) \to M^{\Phi}(\Z^n)$ are bounded. Let us start with the following propositions.
\begin{proposition}\label{pro4}  
Let $\sigma \in M^1(\Z^n \times \T^n)$. Let $\left(\Phi, \Psi \right)$ be a complementary Young pair, and $\Phi$ satisfies a local $\Delta_2$-condition and continuous. Let $g_1 \in  M^{\Psi}(\Z^n)$ and $g_2 \in M^{\Phi}(\Z^n)$. Then, the operator $\mathfrak{L}^{g_1, g_2}_{\sigma}: M^{\Phi}(\Z^n) \to M^{\Phi}(\Z^n)$ is a bounded linear operator, and 
\[\Vert \mathfrak{L}^{g_1, g_2}_{\sigma} \Vert_{\mathcal{B}(M^{\Phi}(\Z^n))} \leq  \Vert \sigma  \Vert_{M^1(\Z^n \times \T^n)} \; \Vert g_1 \Vert_{ M^{\Psi}(\Z^n)} \; \Vert g_2 \Vert_{M^{\Phi}(\Z^n)}. \] 
\end{proposition}
\begin{proof}
Let $f \in M^{\Phi}(\Z^n)$ and $g \in M^{\Psi}(\Z^n)$. Since $(M^{\Phi}(\Z^n))^* \cong M^{\Psi}(\Z^n)$, applying H\"older's inequality, we get
\begin{equation}\label{eq52}
\left| V_{g} f (m, w) \right| \leq \Vert f  \Vert_{M^{\Phi}(\Z^n)} \; \Vert g  \Vert_{M^{\Psi}(\Z^n)}. 
\end{equation} 
For any $f \in M^{\Phi}(\Z^n)$ and $h \in M^{\Psi}(\Z^n)$, applying relations (\ref{eq42}) and (\ref{eq52}), we have
\begin{eqnarray*}
&& \left| \left\langle \mathfrak{L}^{g_1, g_2}_{\sigma}f, h \right\rangle \right| \\
& \leq & \sum_{m \in \Z^n}  \int_{\T^n} |\sigma(m, w)| \; \left|V_{g_1}f(m, w) \right| \left|   V_{g_2}h(m, w) \right| \; dw \\
&& \leq \Vert \sigma \Vert_{M^1(\Z^n \times \T^n)} \; \Vert f \Vert_{M^{\Phi}(\Z^n)} \;  \Vert g_1 \Vert_{M^{\Psi}(\Z^n)} \; \Vert h \Vert_{M^{\Psi}(\Z^n)} \; \Vert g_2 \Vert_{M^{\Phi}(\Z^n)}. 
\end{eqnarray*}
Therefore,
\[\Vert \mathfrak{L}^{g_1, g_2}_{\sigma} \Vert_{\mathcal{B}(M^{\Phi}(\Z^n))} \leq  \Vert \sigma  \Vert_{M^1(\Z^n \times \T^n)} \; \Vert g_1 \Vert_{ M^{\Psi}(\Z^n)} \; \Vert g_2 \Vert_{M^{\Phi}(\Z^n)}. \] 
\end{proof}
Next, using the Schur technique we get an $M^{\Phi}(\Z^n)$-boundedness of the operator $\mathfrak{L}^{g_1, g_2}_{\sigma}$. We obtain a new estimate for $\Vert \mathfrak{L}^{g_1, g_2}_{\sigma} \Vert_{\mathcal{B}(M^{\Phi}(\Z^n))}$. 
\begin{proposition}\label{pro5}
Let $\sigma \in M^1(\Z^n \times \T^n)$ and $g_1, g_2 \in  M^1(\Z^n) \cap \ell^\infty(\Z^n)$. Then, there exists a bounded linear operator $\mathfrak{L}^{g_1, g_2}_{\sigma}: M^{\Phi}(\Z^n) \to  M^{\Phi}(\Z^n) $ for which
\begin{eqnarray*}
&& \Vert \mathfrak{L}^{g_1, g_2}_{\sigma} \Vert_{\mathcal{B}(M^{\Phi}(\Z^n))} \\
& \leq & \max(\Vert g_1 \Vert_{M^1(\Z^n)} \Vert g_2 \Vert_{\ell^\infty(\Z^n)}, \Vert g_1 \Vert_{\ell^\infty(\Z^n)} \Vert g_2 \Vert_{M^1(\Z^n)}) \; \Vert \sigma \Vert_{M^1(\Z^n \times \T^n)}.
\end{eqnarray*}
\end{proposition}
\begin{proof}
We define the function $\mathcal{K}$ on $\Z^n \times \Z^n$ by
\begin{equation}\label{eq53}
\mathcal{K}(k, l) = \sum_{m \in \Z^n}  \int_{\T^n} \sigma(m, w)\; \overline{M_w T_m g_1(l)} \; M_w T_m g_2(k) \; dw.
\end{equation}
Then, we have
\begin{align*}
\mathfrak{L}^{g_1, g_2}_{\sigma}f(k)
& = \sum_{m \in \Z^n}  \int_{\T^n} \sigma(m, w) \; V_{g_1}f(m,w) \; M_w T_m g_2(k) \; dw \\
& = \sum_{m \in \Z^n}  \int_{\T^n} \sigma(m, w) \; \left\langle f, M_w T_m {g_1} \right\rangle \; M_w T_m g_2(k) \; dw \\
& = \sum_{m \in \Z^n}  \int_{\T^n} \sigma(m, w) \; \sum_{l \in \Z^n} f(l) \overline{M_w T_m {g_1} (l)} \; M_w T_m g_2(k) \; dw \\
& = \sum_{l \in \Z^n} f(l) \sum_{m \in \Z^n}  \int_{\T^n}  \sigma(m, w) \; \overline{M_w T_m {g_1} (l)}\; M_w T_m g_2(k) \; dw \\
& = \sum_{l \in \Z^n} \mathcal{K}(k,l) \; f(l).
\end{align*}
For every $l \in \Z^n$, we get
\begin{eqnarray}\label{eq005}
\sum_{k \in \Z^n} | \mathcal{K}(k, l) | 
&\leq & \sum_{k \in \Z^n} \sum_{m \in \Z^n}  \int_{\T^n} |\sigma(m, w)| \left|\overline{ M_w T_m g_1(l)} \right| \; \left|M_w T_m g_2(k) \right| \; dw \nonumber \\
&& \leq \Vert g_1 \Vert_{\ell^\infty(\Z^n)} \; \Vert g_2 \Vert_{M^1(\Z^n)} \; \Vert \sigma \Vert_{M^1(\Z^n \times \T^n)}, 
\end{eqnarray} 
and for every $k \in \Z^n$, we have
\begin{equation}\label{eq006}
\sum_{l \in \Z^n} | \mathcal{K}(k, l) | \leq \Vert g_1 \Vert_{M^1(\Z^n)} \; \Vert g_2 \Vert_{\ell^\infty(\Z^n)} \; \Vert \sigma \Vert_{M^1(\Z^n \times \T^n)}.
\end{equation}
The kernel function $\mathcal{K}(k, l)$ of the localization operator $\mathfrak{L}^{g_1, g_2}_{\sigma}$ satisfies \eqref{eq005} and \eqref{eq006}. Hence, by Schur's lemma (see \cite{fol95}), the operator $\mathfrak{L}^{g_1, g_2}_{\sigma}$ extends to a bounded linear operator $\mathfrak{L}^{g_1, g_2}_{\sigma}: M^{\Phi}(\Z^n) \to  M^{\Phi}(\Z^n)$ with the operator norm 
\begin{eqnarray*}
&& \Vert \mathfrak{L}^{g_1, g_2}_{\sigma} \Vert_{\mathcal{B}(M^{\Phi}(\Z^n))} \\
& \leq & \max(\Vert g_1 \Vert_{M^1(\Z^n)} \Vert g_2 \Vert_{\ell^\infty(\Z^n)}, \Vert g_1 \Vert_{\ell^\infty(\Z^n)} \Vert g_2 \Vert_{M^1(\Z^n)}) \; \Vert \sigma \Vert_{M^1(\Z^n \times \T^n)}.
\end{eqnarray*}
\end{proof} 
\begin{remark}
Based on Proposition \ref{pro5}, it is determined that the bounded linear operator on $M^{\Phi}(\Z^n)$ identified in Proposition \ref{pro4}, is in fact the discrete integral operator on  $M^{\Phi}(\Z^n)$ characterized by the kernel $\mathcal{K}$ as provided in (\ref{eq53}).
\end{remark}

\begin{theorem}\label{th3}
Let $\sigma \in L^1(\Z^n \times \T^n)$. Let $\left(\Phi, \Psi \right)$ be a complementary Young pair, and $\Phi$ satisfies a local $\Delta_2$-condition and continuous. Let $g_1 \in M^{\Psi}(\Z^n)$ and $g_2 \in  M^{\Phi}(\Z^n)$. Then, the operator $ \mathfrak{L}^{g_1, g_2}_{\sigma} \in \mathcal{B}(M^{\Phi}(\Z^n))$, and 
\[ \Vert \mathfrak{L}^{g_1, g_2}_{\sigma} \Vert_{\mathcal{B}(M^{\Phi}(\Z^n))} \leq \Vert \sigma \Vert_{L^1(\Z^n \times \T^n)} \; \Vert g_1 \Vert_{M^{\Psi}(\Z^n)} \; \Vert g_2 \Vert_{M^{\Phi}(\Z^n)}. \]
\end{theorem}
\begin{proof}
Let $f \in M^{\Phi}(\Z^n)$ and $h \in M^{\Psi}(\Z^n)$. Using the duality between $M^{\Phi}(\Z^n)$ and $M^{\Psi}(\Z^n)$, we get 
\begin{eqnarray*}
&& \left| \left\langle \mathfrak{L}^{g_1, g_2}_{\sigma}f, h \right\rangle \right| \\
& \leq & \sum_{m \in \Z^n}  \int_{\T^n} |\sigma(m, w)| \; \left|V_{g_1}f(m, w) \right| \; \left| V_{g_2}h(m, w) \right| \; dw  \\ 
& \leq & \left\Vert \sigma \right\Vert_{L^1(\Z^n \times \T^n)} \Vert f \Vert_{M^{\Phi}(\Z^n)} \Vert g_1 \Vert_{M^{\Psi}(\Z^n)} \Vert h \Vert_{M^{\Psi}(\Z^n)} \Vert g_2 \Vert_{M^{\Phi}(\Z^n)} .
\end{eqnarray*}
Hence,
\[ \Vert \mathfrak{L}^{g_1, g_2}_{\sigma} \Vert_{\mathcal{B}(M^{\Phi}(\Z^n))} \leq \Vert \sigma \Vert_{L^1(\Z^n \times \T^n)} \; \Vert g_1 \Vert_{M^{\Psi}(\Z^n)} \; \Vert g_2 \Vert_{M^{\Phi}(\Z^n)}. \]
\end{proof}

\begin{theorem}
Let $\left(\Phi_i, \Psi_i \right)$ be complementary Young pairs which satisfy local $\Delta_2$-condition, strictly convex and continuous for $i = 1, 2$. Let $\sigma \in L^{\Phi_1, \Phi_2} (\Z^n \times \T^n)$ and $g_1, g_2 \in \mathcal{S}(\Z^n)$. Also, let there exists a constant $x_0 > 0$ such that $\Psi_i(x) \lesssim \Phi_i(x)$ for all $0 \leq x \leq x_0$. Then the localization operator $ \mathfrak{L}^{g_1, g_2}_{\sigma}$ is in $\mathcal{B}(M^{\Phi_1, \Phi_2}(\Z^n))$, and we have 
\[ \Vert \mathfrak{L}^{g_1, g_2}_{\sigma} \Vert_{\mathcal{B}(M^{\Phi_1, \Phi_2}(\Z^n))} \leq \Vert \sigma \Vert_{L^{\Phi_1, \Phi_2} (\Z^n \times \T^n)} . \]
\end{theorem}
\begin{proof}
Let $f \in M^{\Phi_1, \Phi_2}(\Z^n)$ and $h \in M^{\Psi_1, \Psi_2}(\Z^n)$. From the given condition, $\Phi_i, \Psi_i$ satisfy condition (3) of Theorem \ref{th7}. Hence, we have $M^{\Phi_1, \Phi_2}(\Z^n) \subseteq M^{\Psi_1, \Psi_2}(\Z^n)$. Now, applying H\"older's inequality (\ref{eq2}), we obtain
\begin{align*}
\left| \left\langle \mathfrak{L}^{g_1, g_2}_{\sigma}f, h \right\rangle \right| 
& \leq \sum_{m \in \Z^n}  \int_{\T^n} |\sigma(m, w)| \; \left|V_{g_1}f(m, w) \right| \; \left| V_{g_2}h(m, w) \right| \; dw  \\ 
& \leq \left\Vert \sigma \right\Vert_{L^{\Phi_1, \Phi_2} (\Z^n \times \T^n)} \Vert V_{g_1}f \Vert_{L^{\Psi_1, \Psi_2} (\Z^n \times \T^n)}   \Vert V_{g_2}h \Vert_{L^{\Psi_1, \Psi_2} (\Z^n \times \T^n)} \\
& = \left\Vert \sigma \right\Vert_{L^{\Phi_1, \Phi_2} (\Z^n \times \T^n)} \Vert  f \Vert_{M^{\Psi_1, \Psi_2}(\Z^n)}   \Vert  h \Vert_{M^{\Psi_1, \Psi_2}(\Z^n)} \\
& \leq \left\Vert \sigma \right\Vert_{L^{\Phi_1, \Phi_2} (\Z^n \times \T^n)} \Vert  f \Vert_{M^{\Phi_1, \Phi_2}(\Z^n)}   \Vert  h \Vert_{M^{\Psi_1, \Psi_2}(\Z^n)}. 
\end{align*}
Hence, using the duality between $M^{\Phi_1, \Phi_2}(\Z^n)$ and $M^{\Psi_1, \Psi_2}(\Z^n)$, we get
\[ \Vert \mathfrak{L}^{g_1, g_2}_{\sigma} \Vert_{\mathcal{B}(M^{\Phi_1, \Phi_2}(\Z^n))} \leq \Vert \sigma \Vert_{L^{\Phi_1, \Phi_2} (\Z^n \times \T^n)}. \]
This completes the proof.
\end{proof}

\section*{Acknowledgments}
The authors wish to thank the referees for their valuable comments and suggestions that helped to improve the quality of the paper. The second author is partially supported by the XJTLU Research Development Fund (RDF-23-01-027).

\section*{Conflict of interest}
The authors declare that there is no potential conflict of interest regarding the publication of this article.

\section*{Data Availability}
The authors confirm that the data supporting the findings of this study are available within the article and its supplementary materials.

\section*{ORCID}
{\it Aparajita Dasgupta}  https://orcid.org/0000-0001-7093-8158 

{\it Anirudha Poria}  https://orcid.org/0000-0002-0224-3642

\end{document}